\documentclass[12pt,a4paper,english,reqno]{amsart}
\usepackage[a4paper,footskip = 1.5em]{geometry}
\usepackage{amsmath,amssymb,amsthm,mathtools,mathrsfs}
\usepackage{adjustbox,tikz,calc,graphics,babel}
\usepackage{csquotes,enumerate,verbatim}
\usepackage[final]{microtype}
\usepackage[numbers]{natbib}
\usetikzlibrary{shapes.misc,calc,intersections,patterns,decorations.pathreplacing}
\usepackage{hyperref}
\hypersetup{colorlinks=true, linkcolor=blue, citecolor=blue}

\theoremstyle{plain}
\newtheorem*{theorem*}{Theorem}
\newtheorem{theorem}{Theorem}[section]
\newtheorem{lemma}[theorem]{Lemma}
\newtheorem{claim}[theorem]{Claim}

\newtheorem*{claim*}{Claim}

\newtheorem{conjecture}[theorem]{Conjecture}

\theoremstyle{remark}

\def\P{\mathbb{P}}
\def\E{\mathbb{E}}
\def\C{\mathcal}

\let\emptyset\varnothing

\let\originalleft\left
\let\originalright\right
\renewcommand{\left}{\mathopen{}\mathclose\bgroup\originalleft}
\renewcommand{\right}{\aftergroup\egroup\originalright}

\makeatletter
\def\imod#1{\allowbreak\mkern10mu({\operator@font mod}\,\,#1)}
\makeatother
\begin{document}

\title{Separating path systems}

\author{Victor Falgas-Ravry}
\address{Institutionen f\"or matematik och matematisk statistik, Ume{\aa}  Universitet, 901\thinspace87 Ume{\aa}, Sweden}
\email{victor.falgas-ravry@math.umu.se}

\author{Teeradej Kittipassorn}
\address{Department of Mathematical Sciences, University of Memphis, Memphis TN 38152, USA}
\email{t.kittipassorn@memphis.edu}

\author{D\'aniel Kor\'andi}
\address{Departement Mathematik, ETH Z\"urich, R\"amistrasse 101, 8092 Z\"urich, Switzerland}
\email{daniel.korandi@math.ethz.ch}

\author{Shoham Letzter}
\address{Department of Pure Mathematics and Mathematical Statistics, University of Cambridge, Wilberforce Road, Cambridge CB3\thinspace0WB, UK}
\email{s.letzter@dpmms.cam.ac.uk}

\author{Bhargav Narayanan}
\address{Department of Pure Mathematics and Mathematical Statistics, University of Cambridge, Wilberforce Road, Cambridge CB3\thinspace0WB, UK}
\email{b.p.narayanan@dpmms.cam.ac.uk}

\date{12 November 2013}

\begin{abstract}
We study separating systems of the edges of a graph where each member of the separating system is a path. We conjecture that every $n$-vertex graph admits a separating path system of size $O(n)$ and prove this in certain interesting special cases. In particular, we establish this conjecture for random graphs and graphs with linear minimum degree. We also obtain tight bounds on the size of a minimal separating path system in the case of trees.
\end{abstract}

\maketitle

\section{Introduction}
Given a set $S$, we say that a family $\C{F}$ of subsets of $S$ \emph{separates} a pair of distinct elements $x,y \in S$ if there exists a set $A\in \C{F}$ which contains exactly one of $x$ and $y$. If $\C{F}$ separates all pairs of distinct elements of $S$, we say that $\C{F}$ is a \emph{separating system} of $S$.

The study of separating systems was initiated by R\'enyi~\citep{renyi} in 1961. It is essentially trivial that the minimal size of a separating system of an $n$-element set is $\lceil \log_2{n} \rceil$. However, the question of finding the minimal size of a separating system becomes much more interesting when one imposes restrictions on the elements of the separating system. For example, separating systems with restrictions on the cardinalities of their members have been studied by Katona~\citep{katona}, Wegener~\citep{wegener}, Ramsay and Roberts~\citep{ramrob} and K\"undgen, Mubayi and Tetali~\citep{mubtet}, amongst others. Stronger notions of separation as well as other extremal questions about separating systems have also been studied; see, for example, the papers of Spencer~\citep{spencer}, Hansel~\citep{hansel}, and Bollob\'as and Scott~\cite{belascott}.

Another interesting direction involves imposing a graph structure on the underlying ground set and imposing graph theoretic restrictions on the separating family (see, for instance,~\citep{cheng, belascott2}). In this paper, we investigate the question of separating the edges of a graph using paths. Given a graph $G=(V,E)$, we say that a family $\C{P}$ of subsets of the edge set $E(G)$ is a \emph{separating path system of $G$} if $\C{P}$ separates $E(G)$ and every element of $\C{P}$ is a path of $G$. The analogous question of separating the vertices of a graph with paths has also been studied; we refer the reader to~\citep{vertexsep} for details.

Separating path systems arise naturally in the context of network design (see, for instance,~\citep{network1, network2, ladder}). We are presented with a communication network with one (and at most one) defective link and our goal is to identify this link. Of course, one could  test every link, but this is not very efficient; can we do better? A natural test to perform is to send a message between a pair of nodes along a predetermined path; if the message does not reach its intended destination, we conclude that the defective link lies on this path. If we model the communication network as a graph, a fixed set of such tests succeeds in locating any defective link if and only if the corresponding family of paths is a separating path system of the underlying graph. We are naturally led to the  following question: what is the size of a minimal separating path system of a given graph?

\begin{figure}\begin{center}
\trimbox{0cm 0cm 0cm -1.0cm}{ 
\begin{tikzpicture}
\foreach \x in {1,2,3,4,5,6,7,8,9,10,11}
 \node (\x) at (\x, 0) [inner sep=0.5mm, circle, fill=black!100] {};
\draw (1,0) -- (11,0);
\draw (2,1) -- (4,1);
\draw  (9, 1) -- (11,1);
\draw (3,2) -- (6,2);
\draw (7,2) -- (10,2);
\draw (5,3) -- (8,3);
\end{tikzpicture}
}
\end{center}
\caption{A path on $11$ vertices and a separating path system with $5$ paths.}\label{sepsys13-path}

\end{figure}
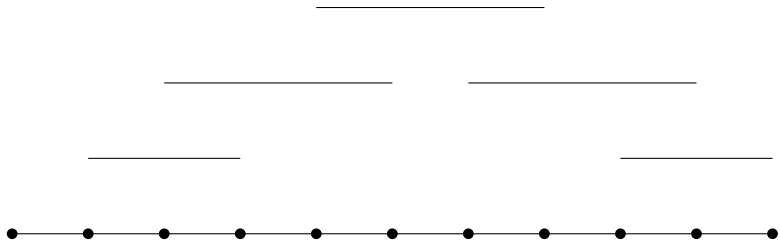

\begin{figure}
\begin{center}
\trimbox{0cm -0.5cm 0cm -0.5cm}{ 
\begin{tikzpicture}
\foreach \x in {1,2,3,4,5,6,7,8,9,10,11}
 \node (\x) at (\x, 5/2) [inner sep=0.5mm, circle, fill=black!100] {};
\foreach \x in {1,2,3,4,5,6,7,8,9,10,11}
 \node (\x) at (\x, 3.5) [inner sep=0.5mm, circle, fill=black!100] {};
 \foreach \x in {1,2,3,4,5,6,7,8,9,10,11}
  \draw (\x, 5/2) -- (\x, 3.5);
\draw (1,2.5) -- (11,2.5);
\draw (1,3.5) -- (11,3.5);

\foreach \x in {1,2,3,4,5,6,7,8,9,10,11}
 \node (\x) at (\x, 0) [inner sep=0.5mm, circle, fill=black!100] {};
\foreach \x in {1,2,3,4,5,6,7,8,9,10,11}
 \node (\x) at (\x, 1) [inner sep=0.5mm, circle, fill=black!100] {};

\draw (1,0) -- (4,0);
\draw (6,0) -- (9,0);
\draw (10,0) -- (11,0);
\draw (4,1) -- (6,1);
\draw (9,1) -- (10,1);
\draw (4,0) -- (4,1);
\draw (6,0) -- (6,1);
\draw (9,0) -- (9,1);
\draw (10,0) -- (10,1);
\end{tikzpicture}
}
\end{center}
\caption{The graph $L_{11}$ and the path $P_A$ corresponding to the subset $A = \{ 4,5,9\}$.}\label{sepsys13-ladder}
\end{figure}
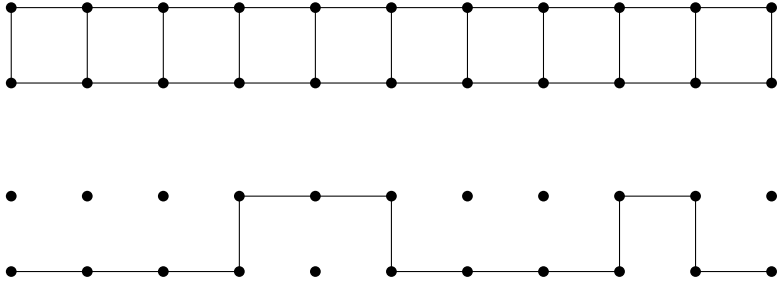

For a graph $G$, let $f(G)$ be the size of a minimal separating path system of $G$. As a separating path system of $G$ is also a separating system of $E(G)$, it follows that $f(G) \ge \lceil \log_2{|E(G)|} \rceil$. In particular, for any connected $n$-vertex graph $G$, $f(G) = \Omega(\log{n})$. With a little work, we can construct graphs that come close to matching this bound. Let $L_n$ be the \emph{ladder} of order $2n$, that is, the Cartesian product of a path of length $n-1$ with a single edge. Given any subset $A$ of $[n-1]$, there is (see Figure~\ref{sepsys13-ladder}) a natural way of mapping $A$ to a path $P_A$ in $L_n$. With this, it is an easy exercise to establish that $f(L_n) = O(\log{n})$; indeed, one can show that $f(L_n) \le 3\log_2{n} + 1$.

A more interesting problem is to determine $f(n)$, the \emph{maximum} of $f(G)$ taken over all $n$-vertex graphs; this question was raised by Gyula Katona in August 2013 at the 5\textsuperscript{th} Eml\'ekt\'abla Workshop in Budapest.

Clearly, at most one edge of a graph can be left uncovered by the paths of a separating path system of the graph; it is thus unsurprising that the question of building small separating path systems is closely related to the well-studied question of covering a graph with paths. It would be remiss not to remind the reader of a beautiful conjecture of Gallai which asserts that every connected graph on $n$ vertices can be decomposed into $\lfloor (n+1)/2 \rfloor$ paths. The following fundamental result of Lov\'asz~\citep{lovasz}, which provides support for Gallai's conjecture, will prove useful; here and elsewhere, by a decomposition of a graph we mean a covering of its edges with edge disjoint subgraphs.

\begin{theorem}
\label{sepsys13-path_decomposition}
Every $n$-vertex graph can be decomposed into at most $n/2$ paths and cycles. Consequently, every $n$-vertex graph can be decomposed into at most $n$ paths. \qed
\end{theorem}

Let $G$ be any graph on $n$ vertices and let $E_1, E_2, \dots, E_k$ be a separating system of the edge set $E(G)$ where $k = \lceil \log_2{|E(G)|} \rceil \le 2\log_2{n}$. Let $G_i$ be the subgraph of $G$ induced by the edges of $E_i$. By Theorem~\ref{sepsys13-path_decomposition}, each $G_i$ may be decomposed into at most $n$ paths. Putting these together, we get a separating path system of $G$ of cardinality at most $kn$. Consequently, we note that $f(n) \le 2n\log_2{n}$.

To bound $f(n)$ from below, let us consider $K_n$, the complete graph on $n$ vertices. Suppose that we have a separating path system $\C{P}$ of $K_n$ with $k$ paths. Note that at most one edge of $K_n$ goes uncovered by the paths of $\C{P}$ and furthermore, at most $k$ edges of $K_n$ belong to exactly one path of $\C{P}$. Since any path of $K_n$ has length at most $n-1$, we deduce that 
\[ k(n-1) \ge 1 + k + 2\left(\binom{n}{2} - k - 1\right)\]
or equivalently, $k \ge n-1-1/n$. Thus, we note that $f(n)\ge n-1$. We believe that the lower bound, rather than the upper bound, is closer to the truth; we make the following conjecture.

\begin{conjecture}
\label{sepsys13-main_conj}
There exists an absolute constant $C$ such that for every graph $G$ on $n$ vertices, $f(G) \le Cn$.
\end{conjecture}

Let us remark that it is not inconceivable that $f(n) = (1+o(1))n$ and Conjecture~\ref{sepsys13-main_conj} is true for every $C>1$. In this paper, we shall prove Conjecture~\ref{sepsys13-main_conj} in certain special cases. Our first result establishes the conjecture for graphs of linear minimum degree.

\begin{theorem}\label{sepsys13-theorem_linear_min_deg}
Let $c>0$ be fixed. Every graph $G$ on $n$ vertices with minimum degree at least $cn$ has a separating path system of cardinality at most $122n/c^2$ for all sufficiently large $n$.
\end{theorem}

Building upon the ideas used to prove Theorem~\ref{sepsys13-theorem_linear_min_deg}, we shall prove Conjecture~\ref{sepsys13-main_conj} for the Erd\H{o}s-R\'enyi random graphs using the fact that these graphs  have good connectivity properties.

\begin{theorem}\label{sepsys13-theorem_random_graphs}
For any probability $p = p(n)$, with high probability, the random graph $G(n,p)$  has a separating path system of size at most $48n$.
\end{theorem}

Note in particular that Theorem~\ref{sepsys13-theorem_random_graphs} implies that Conjecture~\ref{sepsys13-main_conj} is true for almost all $n$-vertex graphs with $C=48$. Using Theorem~\ref{sepsys13-theorem_linear_min_deg}, we shall also establish the conjecture for a class of dense graphs, which includes quasi-random graphs (in the sense of  Chung, Graham and Wilson~\citep{quasi_graph1} and Thomason~\citep{quasi_graph2}) as a subclass.

\begin{theorem}\label{sepsys13-theorem_dense}
Let $c>0$ be fixed and let $G$ be a graph on $n$ vertices such that every subset $U \subseteq V(G)$ of size at least $\sqrt{n}$ spans at least $c|U|^2$ edges. Then $f(G) \le 638n/c^3$ for all sufficiently large $n$.
\end{theorem}

The above results are far from best-possible but we make no attempt to optimise our bounds since it seems unlikely that our methods will yield the best-possible constants. However in the case of trees, we are able to obtain tight bounds.

\begin{theorem}\label{sepsys13-theorem_trees}
Let $T$ be a tree on $n\ge 4$ vertices. Then
\[\left\lceil\frac{n+1}{3}\right\rceil \le f(T) \le \left\lfloor\frac{2(n-1)}{3}\right\rfloor.\]
Furthermore, these bounds are best-possible.
\end{theorem}

We use standard graph theoretic notions and notation and refer the reader to~\citep{belabook1} for terms and notation not defined here. We shall also make use of some well known results about random graphs without proof; see~\citep{belabook2} for details. 

The rest of this paper is organised as follows. In the next section, we describe a general strategy that we adopt to prove Theorems~\ref{sepsys13-theorem_linear_min_deg} and~\ref{sepsys13-theorem_random_graphs}. We then prove Theorems~\ref{sepsys13-theorem_linear_min_deg} and~\ref{sepsys13-theorem_random_graphs} in Sections ~\ref{sepsys13-min_deg} and~\ref{sepsys13-random} respectively. Section~\ref{sepsys13-dense} is devoted to the proof of Theorem~\ref{sepsys13-theorem_dense}. We then prove Theorem~\ref{sepsys13-theorem_trees} in Section~\ref{sepsys13-trees}. We conclude the paper in Section~\ref{sepsys13-conclusion} with a discussion of related questions and problems. For the sake of clarity, we systematically omit floors and ceilings in Sections~\ref{sepsys13-min_deg},~\ref{sepsys13-random} and~\ref{sepsys13-dense}.

\section{A general strategy}\label{sepsys13-strategy}
Theorems~\ref{sepsys13-theorem_linear_min_deg} and~\ref{sepsys13-theorem_random_graphs} are proved similarly, using the following strategy. Let $G_1$ and $G_2$ be subgraphs of $G$ which partition the edge set of $G$. First, we decompose the edges of $G_1$ into at most $3n$ matchings $M_1,\dots,M_{3n}$ as follows. Initially, each $M_i$ is empty; we add the edges of $G_1$ one by one to a suitably chosen matching $M_i$. By Theorem~\ref{sepsys13-path_decomposition}, there exists a path decomposition of $G_1$ into (at most) $n$ paths $P_1,\dots,P_n$. Given an edge $e=xy\in E(G_1)$, let $j$ be such that $e\in P_j$. We add $e$ to a matching $M_i$ which contains no edge of $P_j$ and no edge incident to $x$ or $y$. As the length of $P_j$ is at most $n-1$ and there are at most $2n$ edges incident to either $x$ or $y$, this process is well defined; indeed, we can always find a matching $M_i$ satisfying the required conditions. Note that we have ensured that $|M_i \cap P_j|  \le 1$ for each $1 \le i \le 3n$ and $1\le j \le n$.

Next, for each $1\le i \le 3n$, we find a covering of $E(M_i)$ with paths using edges from $E(G_2)\cup E(M_i)$. These covering paths together with the paths $P_1,\dots,P_n$ separate the edges of $G_1$ from each other and from the edges of $G_2$. To check this, consider an edge $e \in E(G_1)$ such that $e \in P_j$ and $e \in M_i$ and note that $P_j$ separates $e$ from every edge of $G_2$ as well as each edge of $E(G_1) \setminus P_j$, while the path covering $M_i$ separates $e$ from every other edge of $P_j$ since $|M_i \cap P_j|  \le 1$. Repeating this process with the roles of $G_1$ and $G_2$ reversed, we obtain a separating path system of $G$. 

In order to prove the existence of a small separating path system, we shall partition the graph $G$ into $G_1$ and $G_2$ in a way that will enable us to keep the cardinalities of the above coverings small.

\section{Graphs of linear minimum degree}\label{sepsys13-min_deg}
\begin{proof}[Proof of Theorem~\ref{sepsys13-theorem_linear_min_deg}]
Let $G$ be a graph on $n$ vertices with minimum degree at least $cn$, for some fixed $0<c<1$. It is easy to decompose $G$ into two disjoint subgraphs $G_1$ and $G_2$ in such a way that both subgraphs have minimum degree at least $cn/3$. Indeed, one way to do this is to define $G_1$ by randomly selecting each edge of $G$ with probability $1/2$ and to take $G_2$ to be the complement  of $G_1$ in $G$, i.e., $V(G_2) = V(G)$ and $E(G_2) = E(G) \setminus E(G_1)$; the minimum degree conditions follow from the standard estimates for the tail of the binomial distribution.

Following the strategy described in Section~\ref{sepsys13-strategy}, let $P_1,\dots,P_{n}$ be a path decomposition of $G_1$ and let $M_1,\dots,M_{3n}$ be a decomposition of $G_1$ into matchings such that the intersection $M_i\cap P_j$ contains at most one edge for each $i$ and $j$.

Define a graph $H$ on $V(G)$ as follows: two distinct vertices $x,y\in V(G)$ are adjacent in $H$ if they have at least $c^2n/24$ common neighbours in $G_2$. Note that $H$ has no independent set of size $4/c$. Indeed, if $A\subseteq V(G)$ is an independent set in $H$ of size $k=4/c$, then writing $\Gamma(.)$ to denote vertex neighbourhoods, we have 
\begin{align*}
n=|V(G)| &\ge \sum\limits_{x\in A}|\Gamma(x,G_2)|-\sum\limits_{x\neq y\in A}|\Gamma(x,G_2)\cap \Gamma(y,G_2)|\\
&>\frac{kcn}{3}-\frac{k^2c^2n}{48} =(4/3-1/3)n=n
\end{align*}
which is a contradiction.

For each $1\le i\le 3n$, define a sequence of paths in $E(M_i)\cup E(H)$ as follows. Colour the edges of $M_i$  blue and the edges of $H$  red; note that there may be edges coloured both red and blue. Let $Q_{i,1}$ be a longest path alternating between blue and red edges and starting with a blue edge. Having defined $Q_{i,1},\dots, Q_{i,j-1}$, we set \[E_{i,j}=E(M_i)\setminus(E(Q_{i,1})\cup\dots\cup E(Q_{i,j-1})).\] If $E_{i,j} = \emptyset$, we stop. If not, let $Q_{i,j}$ be a longest path alternating between blue edges from $E_{i,j}$ and red edges, starting with a blue edge. Note that we might reuse red edges in this process, but not blue edges.

Since each $Q_{i,j}$ is a longest path, the starting vertices of the paths $Q_{i,j}$ form an independent set in $H$. Thus for each $1\le i\le 3n$, we have at most $4/c$ such paths $Q_{i,j}$ and consequently at most $12n/c$ paths in total. Note that every edge of $G_1$ appears exactly once in one of these $12n/c$ paths as a blue edge. Thus the sum of the lengths of these paths $Q_{i,j}$ is at most $2|E(G_1)|\le n^2$. We split each of the paths $Q_{i,j}$ into paths of length $c^2n/48$, where we allow one of the subpaths to have length less than $c^2n/48$. We thus obtain at most $n^2/(c^2n/48)+12n/c\le 60n/c^2$ red-blue paths. Note that for every red edge $xy$, the vertices $x$ and $y$ have at least $c^2n/24$ common neighbours  in $G_2$. Consequently, we can transform all the red-blue paths into simple paths in $G$: we replace every red edge with a path of length two in $G_2$ with the same endpoints. We can do this because the number of common neighbours in $G_2$ of the ends of a red edge is at least twice the length of the original red-blue path. The family consisting of these paths and the paths $P_1,\dots,P_n$ separates the edges of $G_1$ and has size at most $60n/c^2+n\le 61n/c^2$.

By repeating the above process with the roles of $G_1$ and $G_2$ reversed, we obtain a separating path system of $G$ of size at most $122n/c^2$.
\end{proof}

\section{Random graphs}\label{sepsys13-random}

\begin{proof}[Proof of Theorem~\ref{sepsys13-theorem_random_graphs}]
We use different arguments for different ranges of the edge probability.

\textbf{Case 1: $p \ge 10\log{n}/n$.} Let $G$ be a copy of $G(n,2p)$, where $p\ge 5\log n/n$. We define graphs $G_1$ and $G_2$ on the vertex set of $G$ as follows. We construct $G_1$ by randomly selecting each edge of $G$ with probability $1/2$ and we take $G_2$ to be the complement  of $G_1$ in $G$; clearly, $G_1$ and $G_2$ are edge-disjoint copies of $G(n,p)$.

The following lemma is easily proved using the standard estimates for the tail of the binomial distribution.
\begin{lemma}\label{sepsys13-lem_random_graph_properties}
Let $p\ge 5\log n/n$. Then with high probability, the following assertions hold.
\begin{enumerate}
\item $n^2p/4 \le |E(G(n,p))|\le n^2p$.
\item $G(n,p)$ has minimum degree at least $np/5$.
\item $G(n,p)$ is $np/10$-connected.
\end{enumerate}
\end{lemma}

We will need the notion of a \emph{$k$-linked graph}: a graph is said to be \emph{$k$-linked} if it has at least $2k$ vertices and for every sequence of $2k$ distinct vertices $u_1,\dots,u_k,v_1,\dots,v_k$, there exist vertex disjoint paths $P_1,\dots, P_k$ such that the endpoints of $P_i$ are $u_i$ and $v_i$. Bollob\'as and Thomason ~\cite{linked_graphs} showed that every $22k$-connected graph is $k$-linked. This was later improved by Thomas and Wollan~\cite{linked_graphs_new}, who proved that every $2k$-connected graph on $n$ vertices with at least $5kn$ edges is $k$-linked. From the latter result and Lemma~\ref{sepsys13-lem_random_graph_properties}, we conclude that with high probability, both $G_1$ and $G_2$ are $np/20$-linked.

Following the strategy described in Section~\ref{sepsys13-strategy}, we  find a decomposition of $G_1$ into paths $P_1,\dots,P_n$ and a decomposition of $G_1$ into matchings $M_1,\dots,M_{3n}$ such that the intersection $M_i\cap P_j$ contains at most one edge for each $i$ and $j$.

We decompose each matching $M_i$ into submatchings of size at most $np/20$. Since $G_1$ has at most $n^2p$ edges, we thus obtain at most $23n$ different matchings $M_1',\dots, M_{23n}'$. Now since $G_2$ is $np/20$-linked, we can complete each such matching $M_i'$ into a path using the edges of $G_2$. These paths along with $P_1,\dots,P_n$ constitute a separating path system of $G_1$ of size at most $24n$. Reversing the roles of $G_1$ and $G_2$, we obtain a set of $24n$ paths separating the edges of $G_2$. The union of these two families of paths is a separating path system of $G$ of cardinality at most $48n$.

\textbf{Case 2: $p \le 10/n$.} In this case, with high probability, $G(n,p)$ has at most $20n$ edges, so the edges of $G$ constitute a separating path system of size at most $20n$.

\textbf{Case 3: $10/n \le p \le 10\log{n}/n$.} We begin by collecting together some useful properties of sparse random graphs. First, let us establish some notation: given a graph $G$, write $B_i(v)=B_i(v,G)$ for the set of vertices at (graph-)distance at most $i$ from $v$ and let $\Gamma_i(v) = \Gamma_i(v,G) = B_i(v) \setminus B_{i-1}(v)$. The following lemma is somewhat technical; we defer its proof to the end of the section.

\begin{lemma}\label{sepsys13-lem_sparse_random_graphs_properties}
Let $10\le d\le 10 \log n$. Then with high probability, the following assertions hold for $G=G(n,d/n)$.
\begin{enumerate}
\item \label{sepsys13-item_0}
$G$ has at most $dn$ edges.
\item \label{sepsys13-item_1}
$|\Gamma_i(x)|\le (2d)^i\log n$ for every $x\in V(G)$ and  $i\le n$.
\item \label{sepsys13-item_2}
Every set of $i\le \sqrt{n}$ vertices spans at most $2i$ edges. Furthermore, every set of $i\le 10\log\log n$ vertices spans at most $i$ edges.
\item \label{sepsys13-item_3}
If $G'$ is a subgraph of $G$ with minimum degree at least $10$, then $|\Gamma_{i}(x,G')|\ge 2^i$ for every $x\in V(G')$ and every $1\le i\le 10\log\log n$.
\item \label{sepsys13-item_4}
Let $G'$ be a subgraph of $G$ obtained by deleting at most $20d\log n$ vertices and edges and let $l=3 \log\log n$. For every pair of vertices $x,y\in V(G')$ such that $|B_l(x,G')|,|B_l(y,G')|\ge (\log n)^{3},$ there is a path of length at most $2\log n$ between $x$ and $y$ in $G'$.
\end{enumerate}
\end{lemma}

The \emph{$k$-core} of a graph is its largest induced subgraph with minimum degree at least $k$. Let $H$ be the $15$-core of $G=G(n,p)$ and let $d = np$. By Theorem~\ref{sepsys13-path_decomposition}, we can decompose $H$ into $n$ paths. Since by Lemma~\ref{sepsys13-lem_sparse_random_graphs_properties}(\ref{sepsys13-item_0}) there are at most $dn$ edges of $G$, we can decompose these $n$ paths into at most $2n$ subpaths $Q_1, \dots, Q_{2n}$, each of which has length at most $d$.

Let $l=3\log\log n$. We shall define a collection of at most $2n$ matchings in $H$ of size $d$ each using the paths $Q_1,\dots,Q_{2n}$. Each of these matchings will consist of $d$ edges $e_1, \dots, e_d$ chosen from some $d$ distinct paths $Q_{i_1}, \dots, Q_{i_d}$ which have the additional property that for every $j\neq j'$ and every $x\in V(Q_{i_j})$ and $x'\in V(Q_{i_{j'}})$ we have $B_l(x)\cap B_l(x')=\emptyset$.

We begin with a collection of paths $R_1, \dots, R_{2n}$ which we modify as we go along. Initially we set $R_i=Q_i$ for every $i$. We  define our collection of matchings in $H$ in a sequence of rounds. 

At the beginning of a round, if we have fewer than $2\sqrt{n}$ nonempty paths $R_i$, we stop. Otherwise, we select $d$ of the $R_i$ (in a way we specify below), remove the initial edge from each of these paths and use these $d$ removed edges to form a matching of size $d$. To choose our $d$ paths $R_{i_1}, \dots, R_{i_d}$ we proceed as follows. Let $R_{i_1}$ be any nonempty path. Now, assume that we have chosen $R_{i_1},\dots,R_{i_{t-1}}$, where $t\le d$. Let $N_t=\bigcup_x B_{2l+1}(x)$, where the union is taken over all  $x\in V(Q_{i_1})\cup\dots\cup V(Q_{i_{t-1}})$. From Lemma~\ref{sepsys13-lem_sparse_random_graphs_properties}(\ref{sepsys13-item_1}), we see that 
\[\vert N_t \vert < (t-1) d (2d)^{2l+2}\log n < (2d)^{2l+4} \log n <\sqrt{n}.\]
Thus by Lemma~\ref{sepsys13-lem_sparse_random_graphs_properties}(\ref{sepsys13-item_2}), $N_t$ spans at most $2\sqrt{n}$ edges. Since we started the round with more than $2\sqrt{n}$ nonempty paths, there is a path which contains no edge induced by the  vertex set $N_t$; let $R_{i_t}$ be any such a path. We repeat the procedure until the $d$ paths $R_{i_1}, R_{i_2}, \dots, R_{i_d}$ have been obtained. Clearly, the matchings defined by this process are disjoint and of size $d$, so there are at most $n$ of them; denote them by $M_1,\dots,M_{n}$.

In Lemma~\ref{sepsys13-lem_main_sparse_random_graphs} (stated below), we show that for each such matching $M_i$, there is a path containing $E(M_i)$ and avoiding, for every $e\in E(M_i)$, the other edges of the path $Q\in \{Q_1, Q_2\dots, Q_{2n}\}$ containing $e$.

We leave it to the reader to verify that we then obtain a separating system of size at most $19n$ by taking the union of
\begin{enumerate}
\item
the edges $E(G)\setminus E(H)$ of which there are at most $15n$,
\item
the paths $Q_1,\dots,Q_{2n}$,
\item
the edges of $H$ which are not covered by the matchings $M_1,\dots,M_{n}$ of which there are at most $2d\sqrt{n}\le n$, and
\item
the set of $n$ paths promised by Lemma~\ref{sepsys13-lem_main_sparse_random_graphs}.
\end{enumerate} 

We now state and prove Lemma~\ref{sepsys13-lem_main_sparse_random_graphs}.
\begin{lemma}\label{sepsys13-lem_main_sparse_random_graphs}
Let $G=G(n,p)$ be a graph satisfying Lemma~\ref{sepsys13-lem_sparse_random_graphs_properties}. Let $S_1,\dots,S_d$ be vertex-disjoint paths of length at most $d$ in the $15$-core $H$ of $G$. Set $l=3\log\log n$ and assume that $B_l(x)\cap B_l(y)=\emptyset$ for every $x\in V(S_i)$ and $y\in V(S_j)$ with $i\neq j$. For each $i$, select an edge $e_i=x_iy_i$ from $S_i$, and set $M=\{e_1, e_2, \dots , e_d\}$. Then there exists a path in $G$ containing all the edges of $M$ and no other edge from $\bigcup_{1\le i\le d} E(S_i)$.
\end{lemma}
\begin{proof}
Write $E' = (\bigcup_{1\le i\le d}E(S_i))\setminus E(M)$, let $G_0$ be the graph on $V(G)$ with edge set $E(G)\setminus E'$, and let $G_1$ be the graph obtained from $G_0$ by deleting $x_1$. Consider the graph $H_1$ on the vertex set $V(H) \setminus \{ x_1 \}$ with edge set $E(H) \cap E(G_1)$. Note that $H_1$ has minimum degree at least $12$, since by removing the vertex-disjoint paths $S_1,\dots,S_d$ and the vertex $x_1$ we decrease vertex degrees in $H$ by at most 3. Thus by Lemma~\ref{sepsys13-lem_sparse_random_graphs_properties}(\ref{sepsys13-item_3}), $|B_l(v,G_1)|\ge (\log n)^{3}$ for every $v\in V(M)$. 

We define vertex-disjoint paths $P_1,\dots,P_{d-1}$ of size at most $2\log n$ as follows. Suppose that we have already defined the paths $P_1, P_2, \dots, P_{i-1}$ for some $i<d$. Set $G_i = G_1\setminus \bigcup_{1\le j < i} V(P_j)$ and let $P_i$ be a shortest path in $G_i$ connecting $y_i$ to a vertex from $\bigcup_{i+1 \le j \le d}\{x_j, y_j\}$. Relabelling the remaining vertices and edges if necessary, assume that this path connects $y_i$ to $x_{i+1}$.

We shall show by induction that $P_i$ has length at most $2\log n$. Assume that we have defined $P_1, \dots, P_{i-1}$. By the inductive hypothesis, note that we may assume that $G_i$ is obtained by removing at most $2d\log n$ vertices  and at most $d^2\le 10d\log n$ edges from $G$. Consequently, the bound on the length of $P_i$ would follow from Lemma~\ref{sepsys13-lem_sparse_random_graphs_properties}(\ref{sepsys13-item_4}) by showing that $|B_l(y_i, G_i)|, |B_l(x_{i+1}, G_i)| \ge (\log n)^{3}$.

First, we claim that $B_l(x_{i+1},G_i) = B_l(x_{i+1},G_1)$. To see this, first note that for every $j \le i-1$, the sets $B_l(y_{j},G)$ and $B_l(x_{i+1},G)$ are disjoint and consequently, so are $B_l(y_{j},G_j)$ and $B_l(x_{i+1},G_j)$. Since $P_j$ is a shortest path from $y_j$ to $\bigcup_{j+1 \le k \le d}\{x_k, y_k\}$, it follows that $V(P_j) \cap B_l(x_{i+1},G_j) = \emptyset$. Hence, $B_l(x_{i+1},G_{j+1}) = B_l(x_{i+1},G_j)$ for every $j \le i-1$. It follows that $|B_l(x_{i+1},G_i)| \ge (\log n)^{3}$, and by the same argument, $ B_l(y_{i},G_{i-1}) = B_l(y_{i},G_1)$. 

Next, by the minimality of $P_{i-1}$, it is clear that the set $V(P_{i-1}) \cap B_l(y_{i},G_{i-1})$ is contained in the set $V_{i-1}'$ of the last $l+1$ vertices of $P_{i-1}$. Let $H_i$ be the subgraph of $H_1$ induced by the vertex subset $V(H_1)\setminus V_{i-1}'$. We deduce from Lemma~\ref{sepsys13-lem_sparse_random_graphs_properties}(\ref{sepsys13-item_2}) that no vertex of $G_1$ has more than two neighbours in $V_{i-1}'$, so $H_i$ has minimum degree at least 10. By Lemma~\ref{sepsys13-lem_sparse_random_graphs_properties}(\ref{sepsys13-item_3}) we then have $|B_{l}(y_i, H_i)|\ge(\log n)^{3}$. Since $V(H_i)\cap B_l(y_i,G_1)\subseteq V(G_i)\cap B_l(y_i,G_1)$, it follows that $B_{l}(y_i, H_i) \subseteq  B_{l}(y_i, G_{i})$. Hence, $|B_{l}(y_i, G_{i})| \ge (\log n)^{3}$ and Lemma~\ref{sepsys13-lem_main_sparse_random_graphs} follows by Lemma~\ref{sepsys13-lem_sparse_random_graphs_properties}(\ref{sepsys13-item_4}).
\end{proof}

We now complete the proof of Theorem~\ref{sepsys13-theorem_random_graphs} by proving Lemma~\ref{sepsys13-lem_sparse_random_graphs_properties}.
\begin{proof}[Proof of Lemma~\ref{sepsys13-lem_sparse_random_graphs_properties}]
Parts~(\ref{sepsys13-item_0}) and~(\ref{sepsys13-item_1}) of Lemma~\ref{sepsys13-lem_sparse_random_graphs_properties} follow easily from the standard Chernoff-type bounds for the tails of binomial random variables. Part~(\ref{sepsys13-item_2}) is also routine: the probability that a given set of $i$ vertices induces $k$ or more edges is at most $\binom{i^2}{k}(d/n)^{k}$ and a straightforward union bound over all sets of $i$ vertices establishes both the claimed statements.

To prove part~(\ref{sepsys13-item_3}), we assume that $G$ satisfies parts~(\ref{sepsys13-item_1}) and~(\ref{sepsys13-item_2}). Let $G'$ be a subgraph of $G$ with minimum degree at least $10$. Let $x\in V(G')$ and write $\Gamma_i=\Gamma_{i}(x,G')$ and $B_i=B_{i}(x,G')$.
\begin{claim}\label{sepsys13-doubling}
$|\Gamma_i|\ge 2|B_{i-1}|$ for $1\le i\le 10\log\log n$.
\end{claim}
\begin{proof}
Note that by part~(\ref{sepsys13-item_1}), we have $|\Gamma_i|\le (2d)^i\log n$ for  $1\le i\le 10\log \log n$, so $|B_i|\le 2(2d)^{i+1}\log n\le  \sqrt{n}.$
So by part~(\ref{sepsys13-item_2}), $B_i$ spans at most $2|B_i|$ edges for every $1\le i\le 10\log \log n$. Since every vertex of $B_{i-1}$ has degree at least $10$ in $G'$ and $B_{i-1}$ spans at most $2|B_{i-1}|$ edges,  there are at least $6|B_{i-1}|$ edges from $B_{i-1}$ to $\Gamma_i$. As $B_{i-1}$ is connected, $B_i$ must span at least $7|B_{i-1}|-1$ edges.
Since $B_i$ spans at most $2|B_i|$ edges, this implies that $|B_i|\ge (7|B_{i-1}|-1)/2\ge 3|B_{i-1}|$, i.e., $|\Gamma_i|\ge 2|B_{i-1}|$.
\end{proof}

Claim~\ref{sepsys13-doubling} implies in particular that $|\Gamma_i|\ge 2^i$ for $i\le 10\log\log n$, proving part~(\ref{sepsys13-item_3}). In order to prove part~(\ref{sepsys13-item_4}), we will need the following.

\begin{claim} \label{sepsys13-claim_big_neighbourhood}
Let $l=3\log \log n$. Let $G'$ be a graph obtained from $G$ by removing at most $20d\log n$ vertices and edges and let $x$ be a vertex of $G'$ satisfying $\vert B_l(x,G')\vert\ge (\log n)^3$. Then with high probability, for every such $G'$ and $x$, there exists an $i < \log n$ such that $|\Gamma_{i}(x,G')|\ge n/2d$.
\end{claim}
\begin{proof}
Write $\Gamma_i=\Gamma_i(x,G')$ and $B_i=B_i(x,G')$, and let $V'$ and $E'$ be the set of vertices and edges removed from $G$ to obtain $G'$. Note that the assumption on $x$ implies in particular that there exists a $k\le l$ such that $|\Gamma_{k}|\ge (\log n)^{5/2}$. By part~(\ref{sepsys13-item_1}), with high probability we also have $|B_k|=o(n/d)$.

We show that with high probability, for every $G'$ and $x$ as above and $i\ge k$, either $|\Gamma_{i+1}|\ge (d/2)|\Gamma_i|$ or $|\Gamma_i|\ge (n/2d)$. Note that this would prove Claim~\ref{sepsys13-claim_big_neighbourhood}.

Conditional on $|\Gamma_i|\le {n/2d}$ and on $|\Gamma_{j+1}|\ge d|\Gamma_j|/2$ for $k\le j<i$, we shall bound from above the probability that $|\Gamma_{i+1}|\le d|\Gamma_i|/2$. Write $A_i=V(G')\setminus (\Gamma_i \cup V' \cup V(E'))$, and let $A_i'$ be the set of those vertices of $A_i$ that are adjacent to some vertex of $\Gamma_i$. It follows that since $|\Gamma_i|\le {n/2d}$ and $|V'|,|E'|\le 20d\log n$, $|A_i|\ge 9n/10$ for all sufficiently large $n$. We shall estimate the probability that $|A_i'|\le (d/2)|\Gamma_i|$, conditional on $|\Gamma_i|\ge (\log n)^{5/2}$ and $|A_i|\ge 9n/10$.
The probability that a particular vertex of $A_i$ is adjacent to some vertex of $\Gamma_i$ is 
\begin{equation*}
1-(1-d/n)^{|\Gamma_i|}\ge \frac{d|\Gamma_i|}{n}-\frac{1}{2}\left(\frac{d|\Gamma_i|}{n}\right)^2\ge\frac{3d|\Gamma_i|}{4n}
\end{equation*}
since no edge adjacent to any vertex of $A_i$ is deleted in passing to $G'$ from $G$. Thus the expected size of $A_i'$ is at least ${27d|\Gamma_i|/40}$. By appealing to the standard bounds for the tail of a binomial random variable, we see that
\begin{align*} 
\P\big[\,|A_i'|\le \frac{d}{2}|\Gamma_i|\,\big] &\le \P\Big[\,|A_i'|\le \frac{3}{4}\E|A_i'|] \,\Big]\\
&\le \exp{(-\E|A_i'|/32)}\\
&\le \exp{(-d(\log n)^{5/2}/100)}.
\end{align*}
Since we have $2^{O(d(\log n)^2)}$ choices for $G'$, $x$ and $i$, this implies that $\vert A_i'\vert \ge (d/2) \vert \Gamma_i \vert$ with high probability, as required.
\end{proof}
We now complete the proof of part~(\ref{sepsys13-item_4}) of Lemma~\ref{sepsys13-lem_sparse_random_graphs_properties}. Using Claim~\ref{sepsys13-claim_big_neighbourhood}, we can find $s,t < \log n$ such that $|\Gamma_s(x,G')|, |\Gamma_t(y,G')|\ge n/2d$.
If $B_s(x,G')\cap B_t(y,G')\neq\emptyset$, the assertion of part~(\ref{sepsys13-item_4}) follows. Otherwise, note that the probability that there are no edges between $\Gamma_s(x,G')$ and $\Gamma_t(y,G')$ is at most $(1-d/n)^{n^2/4d^2 - 20d\log{n}}\le e^{-n/5d}$. Since we have $2^{O(d(\log n)^2)}$ choices for $x$, $y$ and $G'$, this implies that with high probability, the assertion of part~(\ref{sepsys13-item_4}) holds.
\end{proof}

We have established Lemma~\ref{sepsys13-lem_sparse_random_graphs_properties}, thus completing the proof of Theorem~\ref{sepsys13-theorem_random_graphs}.
\end{proof}

\section{Dense graphs}\label{sepsys13-dense}
\begin{proof}[Proof of Theorem~\ref{sepsys13-theorem_dense}]
Let $c>0$ and let $G$ be a graph on $n$ vertices such that for every $k\ge \sqrt{n}$, every set of $k$ vertices spans at least $ck^2$ edges.

We define a sequence of subgraphs $G=G_0\supseteq G_1\supseteq \dots \supseteq G_{l-1}$ and a related sequence of graphs $H_1,H_2,\dots,H_l$ as follows. Start by setting $G_0=G$. If $|V(G_{i-1})|\le\sqrt{n}$, we stop and take $H_i = G_{i-1}$. Otherwise, we take $H_i$ to be the $(c|G_{i-1}|/2)$-core of $G_{i-1}$ and define $G_i$ to be the graph induced by $V(G_{i-1})\setminus V(H_i)$. Note that the sets $V(H_i)$ form a  partition of $V(G)$.

Let us write $g_i$ and $h_i$ respectively for the number of vertices of $G_i$ and $H_i$. It is well known that the $k$-core of a graph can be found by removing vertices of degree at most $k-1$, in arbitrary order, until no such vertices exist. So the number of edges removed from $G_{i-1}$ to obtain its $(cg_{i-1}/2)$-core is at most $cg_{i-1}^2/2$. Thus, at least $cg_{i-1}^2/2$ edges remain, so $h_i\ge \sqrt{c}g_{i-1} \ge cg_{i-1}\ge cg_i$.

We first separate the internal edges of the graphs $H_i$. Note that $H_i$ has minimum degree at least $ch_i/2$. So we conclude from Theorem~\ref{sepsys13-theorem_linear_min_deg} that $H_i$ has a separating system of size at most $488h_i/c^2$ for every $1\le i<l$. Also, since $V(H_l)\le \sqrt{n}$, we may separate the edges of $H_l$ trivially (by adding each edge individually to our separating path system); this contributes at most $n$ paths. Since the graphs $H_i$ are  pairwise vertex disjoint, we may separate the internal edges of the $H_i$ using at most $488n/c^2$ paths.

It remains to separate the crossing edges between the $H_i$. For $1\le i<l$, let $E_i$ be the set of edges of the form $xy$ where $x\in V(H_i)$ and $y\in V(G_i)$, and let $E_i'$ be the set of such edges $xy$ where $y$ has at least 3 neighbours in $H_i$. Note that every edge of $G$ not contained in any of the $H_i$ is contained in one of the $E_i$. 

We define a $g_i$-edge-coloured multigraph $F_i$ on the vertex set of $H_i$ as follows. If $v\in G_i$ has at least three neighbours in $H_i$, say $x_{1},\dots,x_k$, we add the edges $x_1x_2,x_2x_3,\dots,x_kx_1$ to $F_i$ and colour these edges with the colour $v$; in other words, we add a $v$-coloured cycle through the neighbours of $v$. Note that the degree of every vertex of $F_i$ (as a multigraph) is at most $2g_i$ and every colour class contains at most $h_i$ edges. Since each edge has at most $4g_i + h_i\le 5h_i/c$ edges which are either incident to it or from the same colour class, we can, as in Section~\ref{sepsys13-strategy}, decompose $F_i$ into at most $5h_i/c$ rainbow matchings $M_1,\dots,M_{5h_i/c}$, where by a rainbow matching, we mean a matching containing at most one edge from each colour class.

We now construct another sequence of rainbow matchings decomposing $F_i$ with the following property. Denote by $e_1,\dots,e_m$ the edges of $M_j$ and let $v_1,\dots,v_m$ be their respective colour classes. Let $\alpha_k$ and $\beta_k$ be the two neighbours of $e_k$ in the cycle whose edges have colour $v_k$. In our second sequence of rainbow matchings, we would like the  matching containing $e_k$ to avoid $e_t$, $\alpha_t$ and $\beta_t$ for every $t\neq k$. Since each edge has to avoid at most $4g_i + h_i + 3h_i\le 8h_i/c$ other edges, we can find such a decomposition into at most $8h_i/c$ matchings, say, $M_{5h_i/c+1}, \dots, M_{13h_i/c}$.

We now mimic the proof of Theorem~\ref{sepsys13-theorem_linear_min_deg}. Let us define a graph $H_i'$ on $V(H_i)$ where we join two vertices if they have more than $c^2h_i/24$ common neighbours in $H_i$. For each $1 \le j \le 13h_i/c$, we can find a collection of at most $4/c$ paths whose edges alternate between those of $M_j$ and $H_i'$ which cover each edge of $M_j$ exactly once; we obtain $52h_i/c^2$ such paths in total. We divide these paths into subpaths of length at most $c^2h_i/48$ each, resulting in a collection of at most $96h_i/c^3 + 52h_i/c^2$ paths. Each such path can be transformed into a path in $G$ by replacing every edge from $H_i'$ with a suitably chosen path of length two in $H_i$ and every coloured edge $e=xy$ from $M_j$ with the path $xvy$ where $v$ is the colour of $e$. Since the matchings are rainbow matchings, these paths are guaranteed to be simple.  It is easy to see that the collection of $96h_i/c^3 + 52h_i/c^2$ paths defined above separates $E_i'$.

It remains to separate the edges of $E_i\setminus E_i'$ for $1\le i< l$. Note that there are at most $2(g_1+\dots+g_l)\le 2(h_1+\dots+h_l)/c \le 2n/c$ such edges; we add each such edge to our separating path system.

It is easy to check that we have constructed a separating path system of $G$ of cardinality at most $488n/c^2 + 96n/c^3 + 52n/c^2 + 2n/c \le 638n/c^3$. The result follows.
\end{proof}

\section{Trees}\label{sepsys13-trees}

We begin by collecting together a few simple observations into the following lemma.

\begin{lemma}\label{sepsys13-lem_leaves_deg2} 
Let $T$ be a tree on $n\ge 3$ vertices, and let $\C{P}$ be a separating path system of $T$. Then the following assertions hold.
\begin{enumerate}
\item \label{sepsys13-t:1}With the exception of at most one leaf, every leaf of $T$ is an endpoint of a path in $\C{P}$.
\item \label{sepsys13-t:2}If a path in $\C{P}$ has two leaves $u$ and $v$ as its endpoints, then there must be at least one path in $\C{P}$ which has exactly one of $u$ and $v$ as an endpoint.
\item \label{sepsys13-t:3}Every vertex of degree two in $T$ is an endpoint of a path in $\C{P}$.
\end{enumerate}
\end{lemma}
\begin{proof}
Clearly, a leaf must be an endpoint of any path through it. Since $\C{P}$ separates $E(T)$, there is at most one edge of $T$ which is not covered by any path in $\C{P}$. As $n\ge 3$, $T$ does not consist of a single edge and thus at most one leaf of $T$ is visited by no path in $\C{P}$. This establishes part~(\ref{sepsys13-t:1}).

Suppose that we have a path $P\in \C{P}$ having two leaves $u,v \in V(T)$ as its endpoints. Let $e_u$ and $e_v$ be the edges incident to $u$ and $v$ respectively. Since $\C{P}$ separates $E(T)$, there must be some path $P'\in \C{P}$ containing exactly one of $e_u$ and $e_v$. This establishes part~(\ref{sepsys13-t:2}).

Suppose that $v$ is a vertex of degree two in $T$; let $e_1$ and $e_2$ be the two edges of $T$ incident to $v$. Since $\C{P}$ separates $E(T)$, there must be some path $P\in \C{P}$ containing exactly one of $e_1$ and $e_2$. Since $v$ has degree two, it must be an endpoint of this path $P$. This establishes part~(\ref{sepsys13-t:3}). 
\end{proof}

We split the proof of Theorem~\ref{sepsys13-theorem_trees} into two parts.

\begin{proof}[Proof of the lower bound in Theorem~\ref{sepsys13-theorem_trees}]
Let $T$ be a tree on $n\ge 4$ vertices and let $\C{P}$ be a separating system of $T$. We shall show that $|\C{P}| \ge \lceil{(n+1)/3}\rceil$.

First, suppose that there is a leaf $v$ which is not the endpoint of any path in $\C{P}$. Let $e_v$ be the edge of $T$ incident to $v$. Since $\C{P}$ separates $E(T)$, $e_v$ is the unique edge of $T$ not covered by any path in $\C{P}$.  Delete $v$ from $T$ to obtain a tree $T'$ on $n-1\ge 3$ vertices.

The family $\C{P}$ both covers and separates $E(T')$. From Lemma~\ref{sepsys13-lem_leaves_deg2}, we note that every leaf and every vertex of degree two of $T'$ is the endpoint of at least one path from $\C{P}$. Furthermore, we know that if a path from $\C{P}$ has a pair of leaves for its endpoints, then at least one of those leaves is the endpoint of at least one other path from $\C{P}$.

Let $d_1$ and $d_2$ denote the number of leaves and degree two vertices of $T'$. We claim that $\C{P}$ contains at least $(2d_1+d_2)/3$ paths. To see this, start by placing a red token on every leaf and a blue token on every vertex of degree two in $T'$. We then iterate through the paths of $\C{P}$ in some order and in each iteration, we remove whatever tokens there are at the endpoints of the current path. If both the tokens removed are red, then we know that both the endpoints, say $u$ and $v$, of the current path are leaves and that at least one of them, say $u$, is the endpoint of a different path; we then place a blue token on $u$. Writing $R$ and $B$ respectively for the number of red and blue tokens remaining on the tree, we see that the quantity $2R + B$ does not decrease by more than three in any iteration. Since $\C{P}$ is a separating path system, all the tokens must have been removed by the end of the procedure. It follows that
\[ \vert\C{P} \vert \ge \frac{2d_1+d_2}{3}. \]
Now, note that
\[ 2e(T')=2(n-2)\ge d_1 + 2d_2+ 3(n-1-d_1-d_2), \]
which we can rearrange to get
\[ \frac{2d_1+d_2}{3}\ge \frac{n+1}{3}. \] 
Taken together, these inequalities show that $\vert \C{P}\vert \ge (n+1)/3$.

If on the other hand every leaf of $T$ is the endpoint of some path from $\C{P}$, then by repeating the argument above with $T$ instead of $T'$, we find that $\vert \C{P}\vert \ge (n+2)/3$. We know from Lemma~\ref{sepsys13-lem_leaves_deg2} that these are the only two possibilities; consequently, we are done.

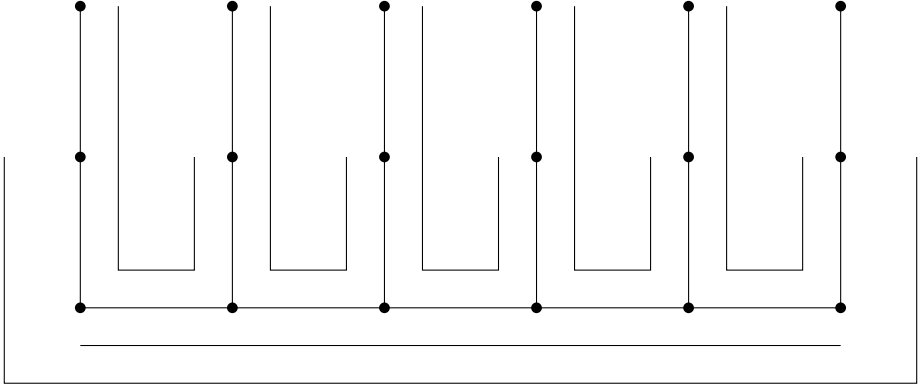
\begin{figure}
\begin{center}
\begin{tikzpicture}[scale=2]
\foreach \x in {1,2,3,4,5,6}
  \node (\x) at (\x, 0) [inner sep=0.5mm, circle, fill=black!100] {};
\foreach \x in {1,2,3,4,5,6}
 \node (1\x) at (\x, 1) [inner sep=0.5mm, circle, fill=black!100] {};
\foreach \x in {1,2,3,4,5,6}
  \node (2\x) at (\x, 2) [inner sep=0.5mm, circle, fill=black!100] {};
\foreach \x in {1,2,3,4,5,6}
  \draw (\x,0) -- (\x,2);
\draw (1,0) -- (6,0);

\draw (0.5, 1) -- (0.5, -0.5) -- (6.5, -0.5) -- (6.5, 1);
\draw (1, -0.25) -- (6, -0.25);

\foreach \x in {1,2,3,4,5}
\draw (\x+0.25, 2) -- (\x+0.25, 0.25) -- (\x+0.75, 0.25) -- (\x+0.75, 1);

\end{tikzpicture}
\end{center}
\caption{A hair-comb of order $18$ and a separating system of $7$ paths.}\label{sepsys13-comb}

\end{figure}

To see that this lower bound is best-possible, consider the family of \emph{hair combs}, where the hair comb of order $3n$ is obtained by starting with a \emph{spine} consisting of a path of length $n-1$ and then attaching a path of length two to each vertex of the spine. It is an easy exercise to show that this lower bound is tight for hair combs. (See Figure~\ref{sepsys13-comb} for an example of an optimal separating path system.)
\end{proof}

We now turn our attention to the second part of the proof of Theorem~\ref{sepsys13-theorem_trees}.
\begin{proof}[Proof of the upper bound in Theorem~\ref{sepsys13-theorem_trees}]
We shall show by induction on $n=\vert V(T)\vert$ that $f(T) \le \lfloor 2(n-1)/3\rfloor$. 

There is, up to isomorphism, only one tree of order $n$ for each of $n=1,2,3$, namely the path of length $n-1$. It is trivial to check that the claim holds for these trees.

Let $T$ be a tree of order $n>3$. If $T$ is a path, then it is easy to show using Lemma~\ref{sepsys13-lem_leaves_deg2} that $f(T) \le \lceil(n-1)/2\rceil$ and it is easy (see Figure~\ref{sepsys13-path}) to construct a separating path system matching this bound. Since $\lceil(n-1)/2\rceil \le \lfloor 2(n-1)/3\rfloor$ for all $n\ge 4$, we may suppose that $T$ is not a path; hence, $T$ must contain at least one vertex $u$ with three distinct neighbours, say $v_1$, $v_2$ and $v_3$. Contract the edges $uv_1$, $uv_2$ and $uv_3$ to obtain a new tree $T'$ on $n-3$ vertices.

We find a separating path system $\C{P}'$ of $T'$ of size at most $2(n-4)/3$. We may think of $\C{P}'$ as a family of paths of $T$ since paths in $T'$ map to paths in $T$ in a natural way: a path in $T'$ is lifted up to a path in $T$ with the same endpoints (where we identify the vertex resulting from the contraction of $u$, $v_1$, $v_2$ and $v_3$ with $u$). Consider the family
\[\C{P}=\C{P}'\cup \{v_1uv_2, v_2uv_3\}.\]
Since $\C{P}'$ separates $E(T')$, it readily follows that $\C{P}'$, when viewed as a family of paths of $T$, separates $E(T)\setminus \{uv_1, uv_2, uv_3\}$. The two paths $v_1uv_2$ and $v_2uv_3$ then separate $uv_1$, $uv_2$ and $uv_3$ from each other and from the rest of $E(T)$. Thus, 
\[ \vert \C{P} \vert \le \frac{2(n-4)}{3}+2=\frac{2(n-1)}{3}. \]
We are done by induction.

To see that this upper bound is best-possible, consider the family of \emph{stars}, where the star of order $n$ consists of a single internal vertex joined to $n-1$ leaves. By mimicking the proof of the lower bound using Lemma~\ref{sepsys13-lem_leaves_deg2}, it is an easy exercise to verify that the upper bound is tight for stars.
\end{proof}

\section{Conclusion}\label{sepsys13-conclusion}
There remain a number of interesting questions which merit investigation. While the main open problem of course is to establish that $f(n) = O(n)$, there are many other attractive related extremal questions. For instance, it would be interesting to determine the value of $f(K_n)$ exactly; one can also ask the same question for the the $d$-dimensional hypercube $Q_d$. It is easy to cover $Q_d$ with $d-1$ ladders, so $f(Q_d)=O(d^2)$. On the other hand, we know from the information theoretic lower bound that $f(Q_d)=\Omega(d)$. It would be interesting to nail down the exact value of $f(Q_d)$.

A different question, though of a similar flavour, raised by Bondy~\citep{Bondy} and answered by Li~\citep{HaoLi}, is that of finding \emph{perfect path double covers}, i.e., a set of paths of a graph such that each edge of the graph belongs to exactly two of the paths and each vertex of the graph is an endpoint of exactly two of the paths. We suspect that the tools developed to tackle this problem and its variants might prove useful in attempting to establish that $f(n) = O(n)$.

\section*{Acknowledgements} Some of the research in this paper was carried out while the authors were participating in the 5\textsuperscript{th} Eml\'ekt\'abla Workshop in Budapest, Hungary.  This research was continued while the fourth and fifth authors were visitors at the University of Memphis. The authors are grateful to the organisers of the Eml\'ekt\'abla Workshop for their hospitality, and the fourth and fifth authors are additionally grateful for the hospitality of the University of Memphis.

\section*{Note added in proof} Shortly after this article was submitted, Balogh, Csaba, Martin and Pluh\'ar~\citep{parallelwork}, working independently, announced some results on a similar problem. Given a graph $G$, they consider the problem of finding a family of paths $\C{P}$ such that for every pair of edges $e,f\in E(G)$, there exist $P_e, P_f \in \C{P}$ such that $e \in P_e, f\notin P_e$ and $f \in P_f, e \notin P_f$. The methods developed in this paper are applicable to their notion of separation in certain cases; it is easy to check that the separating path systems constructed in the proofs of Theorems~\ref{sepsys13-theorem_linear_min_deg},~\ref{sepsys13-theorem_random_graphs} and~\ref{sepsys13-theorem_dense} satisfy their notion of separation as well. However, if we are interested in exact results (as we are in the case of trees and forests, for instance), then it would seem that the two notions of separation lead to different behaviours.

\bibliographystyle{amsplain}
\bibliography{sep_path_sys}

\end{document}